\documentclass[12pt]{article}
\usepackage{latexsym,amssymb,doublespace,amsmath}

\topskip  \newtheorem{theorem}{Theorem}[section]
\newtheorem{proposition}[theorem]{Proposition}
\newtheorem{corollary}[theorem]{Corollary}
\newtheorem{definition}[theorem]{Definition}

\newtheorem{lemma}[theorem]{Lemma}

\newenvironment{proof}{\medskip\noindent{\it Proof.\ }}{\hfill \mbox{$\Box$}\medskip}

\begin{document}
\pagenumbering{arabic}
\def\llim{\lim_{n\rightarrow\infty}}
\def\ls{\leq}
\def\gs{\geq}
\def\LL{\frak L}
\def\qq{{\bold q}}
\def\txx{{\frac1{2\sqrt{x}}}}
\def\B{\Box}
\def\BB{\Box\hspace{-2.5pt}\Box}
\def\md{\mbox{\,{{\footnotesize mod}}\,}}
\title{\sc generalizations of some identities
involving the fibonacci numbers}
\maketitle
\begin{center}Toufik Mansour \footnote{Research financed by EC's
IHRP Programme, within the Research Training Network "Algebraic
Combinatorics in Europe", grant HPRN-CT-2001-00272}
\end{center}
\begin{center}
Department of Mathematics, Chalmers University of Technology,
S-412~96 G\"oteborg, Sweden\\
{\tt toufik@math.chalmers.se}
\end{center}
\section{Introduction}
The generalized Fibonacci and Lucas numbers are defined by
\begin{equation}\label{defuv}
U_n(p,q)=\frac{\alpha^n-\beta^n}{\alpha-\beta},\quad V_n(p,q)=\alpha^n+\beta^n,
\end{equation}
where $\alpha=\frac{1}{2}(p+\sqrt{p^2-4q})$ and
$\beta=\frac{1}{2}(p-\sqrt{p^2-4q})$. Clearly, $U_n(p,q)$ and
$V_n(p,q)$ are the usual Fibonacci and Lucas sequences $\{F_n\}$
and $\{L_n\}$ when $p=1$ and $q=-1$.
\begin{definition}\label{defsum}
Let $d\geq0$. For any $n\geq0$, we define
$$s_d(n;p,q;k)=\sum_{j_1+j_2+\cdots+j_d=n}\,\prod_{i=1}^d U_{k\cdot j_i}(p,q).$$
\end{definition}
For the Fibonacci numbers, Zhang \cite{Z}
found the following identities:
\begin{equation}\label{iz1}
s_2(n;1,-1;1)=\frac{1}{5}((n-1)F_n+2nF_{n-1}),\quad n\geq1,
\end{equation}
\begin{equation}\label{iz2}
s_3(n;1,-1;1)=\frac{1}{50}((5n^2-9n-2)F_{n-1}+(5n^2-3n-2)F_{n-2}),\quad
n\geq2,
\end{equation}
and when $n\geq3$,
\begin{equation}\label{iz3}
s_4(n;1,-1;1)=\frac{1}{150}((4n^3-12n^2-4n+12)F_{n-2}+(3n^3-6n^2-3n+6)
F_{n-3}).
\end{equation}
Recently, Zhao and Wang \cite{ZW} extended these identites to the
case of $\{U_n(p,q)\}$ and $\{V_n(p,q)\}$; for $n\geq1$
\begin{equation}\label{zw1}
s_2(n;p,q;k)=\frac{U_k(p,q)}{V_k^2(p,q)-4q^k}\biggl((n-1)U_{nk}(p,q)V_k(p,q)
-2nq^kU_{(n-1)k}(p,q)\biggr),
\end{equation}
for $n\geq2$,
\begin{equation}\label{zw2}
\begin{array}{ll}
s_3(n;p,q;k)&=\frac{U_k^2(p,q)}{2(V_k^2(p,q)-4q^k)^2}
\biggl((n-1)(n-2)V_k^2(p,q)U_{nk}(p,q)\\[4pt]
&-q^kV_k(p,q)(4n^2-6n-4)U_{(n-1)k}(p,q)\\[2pt]
&+(4n^2-28n+28(n-3)V_k(p,q)+80)U_{(n-2)k}(p,q)\biggr),
\end{array}\end{equation}
and when $n\geq3$,
\begin{equation}\label{zw3}
\begin{array}{ll}
s_4(n;p,q;k)&=\frac{U_k^3(p,q)}{6(V_k^2(p,q)-4q^k)^3}
\biggl(V_k^3(p,q)(n-1)(n-2)(n-3)U_{nk}(p,q)\\[4pt]
&-6q^kV_k^2(p,q)(n-2)(n-3)(n+1)U_{(n-1)k}(p,q)\\[2pt]
&+12q^{2k}V_k(p,q)(n-3)(n^2+n-1)U_{(n-2)k}(p,q)\\[2pt]
&-8q^{3k}n(n^2-4)U_{(n-3)k}(p,q)\biggr).
\end{array}\end{equation}
In this paper, we extend the above conclusions. We establish an identity for
the case $s_d(n;p,q;k)$ for any $d\geq2$.
\section{Main results}
We denote by $G_k(x;p,q)$ the generating function of $\{U_{n\cdot
k}(p,q)\}$, that is, $G_k(x;p,q)=\sum_{n\geq0}U_{n\cdot
k}(p,q)x^n$, where $k$ is a positive integer. Clearly, by
Definition~\ref{defuv} and the geometric formula,
    $$G_k(x;p,q)=\frac{xU_k(p,q)}{1-V_k(p,q)x+q^kx^2}.$$
We define $F_k(x)=F_k(x;p,q)=\frac{G_k(x;p,q)}{x}$. Then
\begin{equation}\label{deff}
F_k(x)=\sum_{n\geq1}U_{n\cdot
k}(p,q)x^{n-1}=\frac{U_k(p,q)}{1-V_k(p,q)x+q^kx^2}.
\end{equation}
\begin{definition}\label{defaa}
Let $a(0,d)=4^d$ for any $d\geq0$, and $a(\ell,0)=0$ for any
$\ell\geq1$. We define $a(\ell,d)$ for $\ell,d\geq1$ by
$a(\ell,d)=4(\ell+1)\cdot a(\ell,d-1)+\ell\cdot a(\ell-1,d-1)$.
\end{definition}
Using this definition we quickly generate the numbers $a(\ell,d)$;
the first few of these numbers are given in Table~1.
$$\begin{array}{llllllll}
d\backslash\ell& 0 & 1 & 2 & 3 & 4 & 5 & 6\\
0 & 1 & 0 & 0 & 0 & 0 & 0 & 0\\
1 & 4 & 1 & 0 & 0 & 0 & 0 & 0\\
2 & 16& 12& 2 & 0 & 0 & 0 & 0\\
3 & 64&112& 48& 6 & 0 & 0 & 0\\
4 &256&960&800&240&24 & 0 & 0\\
5 &1024&7936&11520&6240&1440&120& 0\\
6 &4096&64512&154112&134400&53760&10080&720\\[-20pt]
\end{array}$$
$$\mbox{Table 1: Values of }a(\ell,d)\mbox{ where }0\leq\ell,d\leq6.$$
We can also use Definition~\ref{defaa} to find an explicit formula
for $a(\ell,d)$.
\begin{lemma}\label{lemaa}
For any $\ell,d\geq0$,
$$a(\ell,d)=4^{d-\ell}\sum_{j=0}^\ell(-1)^j\binom{\ell}{j}(\ell+1-j)^d.$$
\end{lemma}
\begin{proof}
By Definition~\ref{defaa} it is easy to see that the lemma holds
for $\ell=0$ or $d=0$. Using induction on $d$ and $\ell$ we get
that
$$\begin{array}{l}
4(\ell+1)\cdot a(\ell,d-1)+\ell\cdot a(\ell-1,d-1)\\
\qquad=(\ell+1)4^{d-\ell}\sum\limits_{j=0}^d(-1)^j\binom{\ell}{j}(\ell+1-j)^{d-1}+
\ell\cdot 4^{d-\ell}\sum\limits_{j=0}^{\ell-1}(-1)^j\binom{\ell-1}{j}(\ell-j)^{d-1}\\
\qquad=4^{d-\ell}\left[(\ell+1)\sum\limits_{j=0}^\ell(-1)^j\binom{\ell}{j}(\ell+1-j)^{d-1}+
\ell\sum\limits_{j=1}^{\ell}(-1)^{j-1}\binom{\ell-1}{j-1}(\ell+1-j)^{d-1}\right]\\[5pt]
\qquad=4^{d-\ell}\left[(\ell+1)^d+\sum\limits_{j=1}^\ell(-1)^j\binom{\ell}{j}(\ell+1-j)^d\right]
=a(\ell,d+1),
\end{array}$$
as requested.
\end{proof}
\begin{definition}\label{defbb}
Let $b(1,d)=(-2)^{d-1}$ for any $d\geq1$, and $b(\ell,1)=0$ for
any $\ell\geq2$. We define $b(\ell,d)$ for $\ell,d\geq2$ by
$b(\ell,d)=b(\ell-1,d-1)-2\ell\cdot b(\ell,d-1)$.
\end{definition}
Using this definition we quickly generate the numbers $b(\ell,d)$;
the first few of these numbers are given in Table~2.
$$\begin{array}{cccccccc}
d\backslash\ell& 1 & 2 & 3 & 4 & 5 & 6 & 7\\
1 & 1 & 0 & 0 & 0 & 0 & 0 & 0\\
2 & -2& 1 & 0 & 0 & 0 & 0 & 0\\
3 & 4 &-6 & 1 & 0 & 0 & 0 & 0\\
4 & -8&28 &-12& 1 & 0 & 0 & 0\\
5 & 16&-120&100&-20&1 & 0 & 0\\
6 &-32&496&-720&260&-30& 1 & 0\\[-20pt]
\end{array}$$
$$\mbox{Table 2: Values of }b(\ell,d)\mbox{ where }0\leq\ell,d\leq6.$$
We can also use Definition~\ref{defbb} to find an explicit formula
for the numbers $b(\ell,d)$.

\begin{lemma}\label{lembb}
For any $\ell,d\geq1$,
$$b(\ell,d)=\frac{(-1)^{d-1}2^{d-\ell}}{(\ell-1)!}
\sum_{j=0}^{\ell-1}(-1)^j\binom{\ell-1}{j}(j+1)^{d-1}.$$
\end{lemma}
\begin{proof}
By Definition~\ref{defbb} it is easy to see that the lemma holds
for $\ell=1$ or $d=1$. Using induction on $d$ and $\ell$ we get
that
$$\begin{array}{l}
b(\ell-1,d-1)-2\ell\cdot b(\ell,d-1)\\
\qquad=\frac{(-1)^{d-2}2^{d-\ell}}{(\ell-2)!}
\sum\limits_{j=0}^{\ell-2}(-1)^j\binom{\ell-2}{j}(j+1)^{d-2}-
\frac{2\ell(-1)^{d-2}2^{d-\ell-1}}{(\ell-1)!}
\sum\limits_{j=0}^{\ell-1}(-1)^j\binom{\ell-1}{j}(j+1)^{d-2}\\
\qquad=\frac{(-1)^{d-1}2^{d-\ell}}{(\ell-1)!}\left[
\ell\sum\limits_{j=0}^{\ell-1}(-1)^j\binom{\ell-1}{j}(j+1)^{d-2}-
(\ell-1)\sum\limits_{j=0}^{\ell-2}(-1)^j\binom{\ell-2}{j}(j+1)^{d-2}
\right]\\[8pt]
\qquad=\frac{(-1)^{d-1}2^{d-\ell}}{(\ell-1)!}\left[
(-1)^{d-1}\ell^{d-1}+\sum\limits_{j=0}^{\ell-2}(-1)^j\left(\ell\binom{\ell-1}{j}-(\ell-1)\binom{\ell-2}{j}\right)(j+1)^{d-2}\right]\\[8pt]
\qquad=\frac{(-1)^{d-1}2^{d-\ell}}{(\ell-1)!}\left[
(-1)^{d-1}\ell^{d-1}+\sum\limits_{j=0}^{\ell-2}(-1)^j\binom{\ell-1}{j}(j+1)^{d-1}\right]\\
\qquad=\frac{(-1)^{d-1}2^{d-\ell}}{(\ell-1)!}\sum\limits_{j=0}^{\ell-1}(-1)^j\binom{\ell-1}{j}(j+1)^{d-1}
=b(\ell,d),
\end{array}$$
as requested.
\end{proof}

Now we introduce a relation that plays the crucial role in the
proof of the main result of this paper.
\begin{proposition}\label{thm1}
Let $d\geq1$. The generating function $F_k(x;p,q)$ satisfies the
following equation:
$$\begin{array}{l}
\sum\limits_{j=0}^d\left[(4q^k)^{d-j}\left(\sum\limits_{i=0}^j(-1)^i\binom{j}{i}(j+1-i)^d\right)\left(\frac{V_k^2(p,q)-4q^k}{U_k(p,q)}\right)^{j}F_k^{j+1}(x;p,q)\right]\\
\quad\quad=\sum\limits_{j=1}^d\left[\frac{(-1)^{d-1}(2q^k)^{d-j}}{(j-1)!}
\left(\sum\limits_{i=0}^{j-1}(-1)^i\binom{j-1}{i}(i+1)^{d-1}\right)\left(V_k(p,q)-2q^kx\right)^jF_k^{(j)}(x;p,q)\right],
\end{array}$$
where $F_k^{(j)}(x;p,q)$ is the $j$th derivative with respect to
$x$ of $F_k(x;p,q)$.
\end{proposition}
\begin{proof}
We define $A=\frac{V_k^2(p,q)-4q^k}{U_k(p,q)}$ and
$B=V_k(p,q)-2q^kx$. Let us prove this theorem by induction on $d$.
Noticing that
\begin{equation}\label{eqonce}
F_k^{(1)}(x;p,q)=\frac{(V_k(p,q)-2q^kx)F_k(x;p,q)}{1-V_k(p,q)x+q^kx^2},
\end{equation}
we get
$$4q^kF_k(x;p,q)+A\cdot F_k^2(x;p,q)=B\cdot F_k^{(1)}(x;p,q),$$
therefore, the theorem holds for $d=1$. Now we suppose that the
theorem holds for $d$, that is,
$$\begin{array}{l}
\sum\limits_{j=0}^da(j,d)q^{(d-j)k}\left(\frac{V_k^2(p,q)-4q^k}{U_k(p,q)}\right)^{j}F_k^{j+1}(x;p,q)\\
\qquad\qquad\qquad\qquad\qquad=\sum\limits_{j=1}^db(j,d)q^{(d-j)k}\left(V_k(p,q)-2q^kx\right)^jF_k^{(j)}(x;p,q).
\end{array}$$
Therefore, derivative this equation with respect to $x$ we have
that
$$\begin{array}{l}
\sum\limits_{j=0}^d(j+1)a(j,d)q^{(d-j)k}A^{j}F_k^{j}(x;p,q)F_k^{(1)}(x;p,q)\\
\qquad\quad=\sum\limits_{j=1}^db(j,d)q^{(d-j)k}B^jF_k^{(j+1)}(x;p,q)-\sum\limits_{j=1}^d2jq^kb(j,d)q^{(d-j)k}B^{j-1}F_k^{(j)}(x;p,q).
\end{array}$$
If multiplying by $B$ and using Equation~\ref{eqonce} then we get
that
$$\begin{array}{l}
\sum\limits_{j=0}^d(j+1)a(j,d)q^{(d-j)k}A^{j+1}F_k^{j+2}(x;p,q)+
\sum\limits_{j=0}^d4(j+1)a(j,d)q^{(d+1-j)k}A^{j}F_k^{j+1}(x;p,q)\\
\qquad\quad=\sum\limits_{j=2}^{d+1}b(j-1,d)q^{(d+1-j)k}B^{j}F_k^{(j)}(x;p,q)-
\sum\limits_{j=1}^{d}2jb(j,d)q^{(d+1-j)k}B^{j}F_k^{(j)}(x;p,q),
\end{array}$$
equivalently,
$$\begin{array}{l}
\sum\limits_{j=0}^{d+1}(ja(j-1,d)+4(j+1)a(j,d))q^{(d+1-j)k}A^{j}F_k^{j+1}(x;p,q)\\
\qquad\qquad\qquad\quad=\sum\limits_{j=1}^{d+1}(b(j-1,d)-2jb(j,d))q^{(d+1-j)k}B^{j}F_k^{(j)}(x;p,q),
\end{array}$$
Therefore, using Definition~\ref{defaa} and Definition~\ref{defbb}
we have that
$$\begin{array}{l}
\sum\limits_{j=0}^{d+1}a(j,d+1)q^{(d+1-j)k}A^{j}F_k^{j+1}(x;p,q)
=\sum\limits_{j=1}^{d+1}b(j,d+1)q^{(d+1-j)k}B^{j}F_k^{(j)}(x;p,q).
\end{array}$$
Hence, using Lemma~\ref{lemaa} and Lemma~\ref{lembb} we get the
desired result.
\end{proof}

By the above proposition, we have the main result of this paper.
\begin{theorem}\label{thm2}
Let $d\geq1$. For any $n\geq d$,
$$\begin{array}{l}
\sum\limits_{j=0}^d\left[(4q^k)^{d-j}\left(\sum\limits_{i=0}^j(-1)^i\binom{j}{i}(j+1-i)^d\right)\left(\frac{V_k^2(p,q)-4q^k}{U_k(p,q)}\right)^{j}s_{j+1}(n+j-d;p,q;k)\right]\\
\,\,=\sum\limits_{j=1}^d\left[\frac{(-1)^{d-1}(2q^k)^{d-j}}{(j-1)!}
\left(\sum\limits_{i=0}^{j-1}(-1)^i\binom{j-1}{i}(i+1)^{d-1}\right)\left(
\sum\limits_{s=0}^jv_{d,j,s}(n)U_{(n+j-d-s)k}(p,q)\binom{j}{s}\right)\right],
\end{array}$$
where $v_{d,j,s}(n)=(-2q^k)^sV_k^{j-s}(p,q)\prod_{i=1}^j(n+j-d-s-i)$.
\end{theorem}
\begin{proof}
If comparing the coefficients of $x^{n-(d+1)}$ on both sides of
Proposition~\ref{thm1} we get the desired result.
\end{proof}

Theorem~\ref{thm2} provides a finite algorithm for finding
$s_d(n;p,q;k)$ in terms of $U_{nk}(p,q)$ and $V_{nk}(p,q)$, since
we have to consider all $s_j(n;p,q;k)$ for $j=1,2,\ldots,d$. The
algorithm has been implemented in Maple, and yields explicit
results for $1\leq d\leq 6$. Below we present several explicit
calculations.
\begin{corollary}{\rm(see Zhao and Wang \cite[Equation~9]{ZW})}\label{cc1}
For any $n\geq1$,
$$s_2(n;p,q;k)=\frac{U_k(p,q)}{V_k^2(p,q)-4q^k}\biggl( (n-1)V_k(p,q)U_{nk}(p,q)
-2nq^kU_{(n-1)k}(p,q)\biggr).$$
\end{corollary}
\begin{proof}
Theorem~\ref{thm2} for $d=2$ yields
$$\begin{array}{l}
4q^ks_1(n-1;p,q;k)+\frac{V_k^2(p,q)-4q^k}{U_k(p,q)}s_2(n;p,q;k)\\
\qquad\qquad\qquad=(n-1)V_k(p,q)U_{nk}(p,q)-2(n-2)q^kU_{(n-1)k}(p,q).
\end{array}$$
Using the fact that $s_1(n;p,q;k)=U_{nk}(p,q)$ we get the desired
result.
\end{proof}

\begin{corollary}{\rm(see Zhao and Wang \cite[Equation~10]{ZW})}\label{cc2}
For any $n\geq2$,
$$\begin{array}{ll}
s_3(n;p,q;k)&=\frac{U_k^2(p,q)}{2(V_k^2(p,q)-4q^k)^2}
\biggl((n-1)(n-2)V_k^2(p,q)U_{nk}(p,q)\\[4pt]
&-2q^k(n-2)(2n+1)V_k(p,q)U_{(n-1)k}(p,q)\\[2pt]
&+4q^{2k}(n-2)(n+2)U_{(n-2)k}(p,q)\biggr).
\end{array}$$
\end{corollary}
\begin{proof}
Theorem~\ref{thm2} for $d=3$ yields
$$\begin{array}{l}
16q^{2k}s_1(n-2;p,q;k)+12q^k\frac{V_k^2(p,q)-4q^k}{U_k(p,q)}s_2(n-1;p,q;k)
+\frac{2(V_k^2(p,q)-4q^k)^2}{U_k^2(p,q)}s_3(n;p,q;k)\\
\qquad=(n-1)(n-2)V_k^2(p,q)U_{nk}(p,q)-2(n-2)(2n-5)q^kV_k(p,q)U_{(n-1)k}(p,q)\\
\qquad+4q^{2k}(n-3)^2U_{(n-2)k}(p,q).
\end{array}$$
Using Corollary~\ref{cc1} with the fact that
$s_1(n;p,q;k)=U_{nk}(p,q)$ we get the desired result.
\end{proof}

Similarly, if applying Theorem~\ref{thm2} for $d$ with using the
formulas of $s_j(n;p,q;k)$ for $j=1,2,\ldots,d-1$, then we get the
following result (in the case $d=4$ see \cite[Equation~11]{ZW}).

\begin{corollary}\label{cc3}
We have

{\rm(i)} For any $n\geq3$,
$$\begin{array}{l}
s_4(n;p,q;k)\\
\quad=\frac{U_k^3(p,q)}{6(V_k^2(p,q)-4q^k)^3}
\biggl(V_k^3(p,q)(n-1)(n-2)(n-3)U_{nk}(p,q)\\[4pt]
\quad-6q^kV_k^2(p,q)(n-2)(n-3)(n+1)U_{(n-1)k}(p,q)\\[2pt]
\quad+12q^{2k}V_k(p,q)(n-3)(n^2+n-1)U_{(n-2)k}(p,q)\\[2pt]
\quad-8q^{3k}n(n^2-4)U_{(n-3)k}(p,q)\biggr).
\end{array}$$

{\rm(ii)} For any $n\geq4$,
$$\begin{array}{l}
s_5(n;p,q;k)\\
\quad=\frac{U_k^4(p,q)}{4!(V_k^2(p,q)-4q^k)^4}
\biggl(V_k^4(p,q)(n-1)(n-2)(n-3)(n-4)U_{nk}(p,q)\\[4pt]
\quad-4q^kV_k^3(p,q)(n-2)(n-3)(n-4)(2n+3)U_{(n-1)k}(p,q)\\[2pt]
\quad+12q^{2k}V_k^2(p,q)(n-3)(n-4)(2n^2+4n-1)U_{(n-2)k}(p,q)\\[2pt]
\quad-8q^{3k}V_k(p,q)(n-4)(2n+1)(2n^2+2n-9)U_{(n-3)k}(p,q)\\[2pt]
\quad+16q^{4k}(n-3)(n-1)(n+1)(n+3)U_{(n-4)k}(p,q)\biggr).
\end{array}$$

{\rm(iii)} For any $n\geq5$,
 $$\begin{array}{l}
s_6(n;p,q;k)\\
\quad=\frac{U_k^5(p,q)}{5!(V_k^2(p,q)-4q^k)^5}
\biggl(V_k^5(p,q)(n-1)(n-2)(n-3)(n-4)(n-5)U_{nk}(p,q)\\[4pt]
\quad-10q^kV_k^4(p,q)(n-2)(n-3)(n-4)(n-5)(n+2)U_{(n-1)k}(p,q)\\[2pt]
\quad+20q^{2k}V_k^3(p,q)(n-3)(n-4)(n-5)(2n^2+6n+1)U_{(n-2)k}(p,q)\\[2pt]
\quad-40q^{3k}V_k^2(p,q)(n-4)(n-5)(n+1)(2n^2+4n-9)U_{(n-3)k}(p,q)\\[2pt]
\quad+80q^{4k}V_k(p,q)(n-5)(n^4+2n^3-10n^2-11n+9)U_{(n-4)k}(p,q)\\[2pt]
\quad-32q^{5k}n(n-4)(n-2)(n+2)(n+4)U_{(n-5)k}(p,q)\biggr).
\end{array}$$
\end{corollary}

From these results, it is very easy to obtain Equations
\ref{iz1}-\ref{iz3}. If $k=1$ and $p=-q=1$, then by using
Corollary~\ref{cc3} together with the recurrence
$F_n=F_{n-1}+F_{n-2}$ we arrive to
$$\begin{array}{l}
\sum\limits_{a+b+c+d+e=n}F_aF_bF_cF_dF_e\\
\qquad\qquad=\frac{1}{4!\cdot5^4}\bigl(3(n-1)(8n^3-5n^2-27n+50)F_n-20n(5n^2-17)F_{n-1}
\bigr)\end{array}$$
$$\begin{array}{l}
\sum\limits_{a+b+c+d+e+f=n}F_aF_bF_cF_dF_eF_f\\
\quad=\frac{1}{5!\cdot 5^4}
\bigl((n-1)(5n^4-70n^3-65n^2+490n+264)F_n+2n(5n^4+5n^2-226)F_{n-1}\bigr).
\end{array}$$

\noindent {\sc 2000 Mathematics Subject Classification}: Primary
11B39, 11B83; Secondary 05A15

\begin{thebibliography}{WWW}
\bibitem{HM}
{\sc A.F.~Horadam and Bro. J. M.~Mahon}, Pell and Pell-Lucas
polynomials, {The Fibonacci Quartely} {\bf 23} (1985) 7--20.
\bibitem{ZW}
{\sc F.~Zhao and T.~Wang}, Generalizations of some identities
involving the Fibonacci numbers, {\em The Fibonacci Quarterly}
{\bf 39} (2001) 165--167.
\bibitem{Z}
{\sc W.~Zhang}, Some identities involving Fibonacci numbers, {\em
The Fibonacci Quarterly} {\bf 35} (1997) 225--229.
\end{thebibliography}
\end{document}